\documentclass[12pt,reqno]{amsart}
\usepackage{amsfonts}
\usepackage{amsmath,amsthm,amssymb,amsfonts,amscd}
\usepackage{mathrsfs}
\usepackage{lipsum}
\usepackage{bbding}
\usepackage{graphicx,latexsym}
\usepackage{hyperref}

\newcommand{\bea}{\begin{eqnarray}}
\newcommand{\eea}{\end{eqnarray}}
\newcommand{\bna}{\begin{eqnarray*}}
\newcommand{\ena}{\end{eqnarray*}}

\numberwithin{equation}{section} 

\renewcommand{\thefootnote}{\fnsymbol{footnote}}
\theoremstyle{plain}

\theoremstyle{definition}

\newcommand\blfootnote[1]{%
  \begingroup
  \renewcommand\thefootnote{}\footnote{#1}%
  \addtocounter{footnote}{-1}%
  \endgroup
}

\renewcommand{\Im}{\operatorname{Im}}

\begin{document}

\title{Shifted convolution sums involving theta series}

\author{Qingfeng Sun}
\address{School of Mathematics and Statistics \\ Shandong University, Weihai
\\ Weihai \\Shandong 264209 \\China}
\email{qfsun@sdu.edu.cn}

\begin{abstract}
Let $f$ be a cuspidal newform (holomorphic or Maass) of arbitrary level and nebentypus and
denote by $\lambda_f(n)$ its $n$-th Hecke eigenvalue.
Let
$$
r(n)=\#\left\{(n_1,n_2)\in \mathbb{Z}^2:n_1^2+n_2^2=n\right\}.
$$
In this paper, we study the shifted convolution sum
$$
\mathcal{S}_h(X)=\sum_{n\leq X}\lambda_f(n+h)r(n), \qquad 1\leq h\leq X,
$$
and establish uniform bounds
with respect to the shift $h$ for $\mathcal{S}_h(X)$.
\end{abstract}


\maketitle

\section{Introduction}
\setcounter{equation}{0}

The shifted convolution sum
\bna
\sum_{n\leq X}a_1(n+h)a_2(n)
\ena
where $a_1(n)$ and $a_2(n)$ are two arithmetic functions and $h\geq 1$ an integer, is an interesting
and important object in analytic number theory and has been studied intensively by many authors
with various applications (see for example \cite{DFI}, \cite{H},
\cite{KMV}, \cite{LS}, \cite{M}).
For the arithmetic function
$$
r_\ell(n)=\#\left\{(n_1,n_2,\ldots,n_\ell)\in \mathbb{Z}^\ell:n_1^2+\cdots+n_\ell^2=n\right\},
$$
which is the $n$-th Fourier coefficient of the modular form $\theta^\ell(z)$, where
$\theta(z)$ is the classical Jacobi theta series
$$
\theta(z)=\sum_{n\in \mathbb{Z}}e(n^2z),
$$
the related shifted convolution problem was first studied by Luo \cite{L}. Precisely,
Luo first established a Voronoi formula for $r_\ell(n)$ and then applying
Poincar\'{e} series reduction and his Voronoi formula proved the following:
Let $\lambda_g(n)$ be the
normalized $n$-th Fourier coefficient of a holomorphic cusp form $g$ of
weight $\kappa$ for $\Gamma_0(N)$.
For $\ell\geq 2$ and $\kappa\geq \frac{\ell}{2}+3$,
\begin{equation}
\sum_{n\leq X}\lambda_g(n+h)r_\ell(n)\ll_{h,g,\ell,\varepsilon}X^{\frac{\ell}{2}-\vartheta_\ell+\varepsilon},
\end{equation}
where
$$
\vartheta_\ell=\frac{\ell-1}{4(\varrho+1)}, \qquad
\varrho=\left\{\begin{array}{ll}\frac{1+\ell}{2},&\mbox{if $\ell$ odd,}\\
1+\frac{\ell}{2},&\mbox{if $\ell$ even.}\end{array}\right.
$$
In particular,
$
\vartheta_2=\frac{1}{12}, \vartheta_3=\frac{1}{6}.
$
Recently, L\"{u}, Wu and Zhai \cite{LWZ} improved Luo's result by the circle method and
showed that (1.1) holds uniformly
for $0\leq h\leq X$ and all $\kappa$, and $\vartheta_\ell$
can be taken as
$
\vartheta_3=\frac{1}{4}, \vartheta_\ell= \frac{1}{2} \,(\ell\geq 4).
$
Moreover, for $N=1$, using some new ideas, they proved that one can take
$
\vartheta_2=\frac{1}{6}, \vartheta_\ell= \frac{2}{3}\, (\ell\geq 6).
$ For other interesting results, see \cite{R}.

\medskip

In this paper, we are concerned with the shifted convolution problem of
theta series with a general $GL(2)$ cusp form, that is (as in \cite{BHM}) a holomorphic
form of integral weight $\kappa\geq 1$, level $D$ and nebentypus $\chi_{D}$, or a
Maass form of weight $0$ or $1$, level $D$, nebentypus $\chi_D$ and Laplace eigenvalue
$1/4+\mu^2$. Suppose that $f$ is a newform and denote by $\lambda_f(n)$
its $n$-th Hecke eigenvalue. Then we have the following result.

\medskip

\noindent {\bf Theorem 1.}\, {\it Let $f$ be a cuspidal newform (holomorphic or Maass) of arbitrary level and nebentypus and
denote by $\lambda_f(n)$ its $n$-th Hecke eigenvalue.
Assume that $\lambda_f(n)\ll n^{\theta+\varepsilon}$ for any $\varepsilon>0$. Then for $2\leq \ell\leq 7$ and $X>1$,
we have
$$
\sum_{n\leq X}\lambda_f(n+h)r_{\ell}(n)
\ll_{f,\ell,\varepsilon}X^{\frac{\ell}{2}-\frac{\ell-1-2\theta}{12}+\varepsilon},
$$
uniformly for $1\leq h\leq X$.
}

The most interesting case of Theorem 1 is $\ell=2$.
Note that the best $\theta$ we known is $\theta=7/64$
(see \cite{K}). Denote $r_2(n):=r(n)$ as usual.
By Theorem 1, we have
\bea
\sum_{n\leq X}\lambda_f(n+h)r(n)
\ll_{f,\varepsilon}X^{1-\frac{25}{384}+\varepsilon}.
\eea

\noindent{\bf Remark 1.}
For $\ell\geq 3$, Theorem 1 is weaker than the results in \cite{LWZ}. In fact,
the argument in Section 3 in \cite{LWZ} yields
$
\vartheta_3=\frac{1}{4}, \,\vartheta_\ell=\frac{1}{2}
$
for $\ell\geq 4$.

To prove Theorem 1, we apply Jutila's variation of the circle method (see \cite{J1})
which gives an approximation for
\[
I_{[0,1]}(x)=\left\{\begin{array}{ll}
1,&x\in [0,1],\\
0,&\mbox{otherwise}.
\end{array}
\right.
\]
Precisely, let $\mathcal{Q}\subset [1,Q]$, $Q>0$ and $Q^{-2}\leq \delta \leq Q^{-1}$. Define
\begin{equation}
\widetilde{I}_{\mathcal{Q},\delta}(x)=\frac{1}{2\delta L}
\sum_{q\in \mathcal{Q}}\sum_{a=1 \atop (a,q)=1}^q
I_{\left[\frac{a}{q}-\delta,\frac{a}{q}+\delta\right]}(x),
\end{equation}
where $L=\sum_{q\in \mathcal{Q}}\phi(q)$.
Then $\widetilde{I}_{\mathcal{Q},\delta}(x)$ is an approximation for $I_{[0,1]}(x)$ in the following sense
(see Lemma 4 in Munshi \cite{M}):
\begin{equation}
\int\limits_0^1\left|1-\widetilde{I}_{\mathcal{Q},\delta}(\beta)\right|^2\mathrm{d}\beta\ll
\frac{Q^{2+\varepsilon}}{\delta L^2}.
\end{equation}
We shall prove Theorem 1 in detail in Section 3.

The advantage of Jutila's version of the circle method is the flexibility of choosing
moduli and this is good for us to study a general $GL(2)$ cusp form.
However, for our problem, there is some lost from (1.4). Thus, for $\ell=2$ and
$f$ a holomorphic cusp form of weight $\kappa$ or a Maass cusp form of
Laplace eigenvalue $\frac{1}{4}+\mu^2$ for $SL(2,\mathbb{Z})$, we further prove the following result.

\medskip

\noindent {\bf Theorem 2.}\, {\it Let $f$ be a holomorphic cusp form of weight $\kappa$ or
a Maass cusp form of
Laplace eigenvalue $\frac{1}{4}+\mu^2$ for $SL(2,\mathbb{Z})$ and
denote by $\lambda_f(n)$ its $n$-th Hecke eigenvalue. Let $\phi(x)$ be a smooth function
which is supported on $[1/2,1]$ and identically equal to 1 on
$[1/2+\Delta^{-1},1-\Delta^{-1}]$ with $\Delta>4$, satisfying
$\phi^{(j)}(x)\ll_j \Delta^{j}$ for all integers $j\geq 0$.
Then for $X>1$ and any $\varepsilon>0$, we have
$$
\sum_{n\geq 1}\lambda_f(n+h)r(n)\phi\left(\frac{n}{X}\right)
\ll_{f,\varepsilon}X^{1-\frac{1}{4}+\varepsilon}\Delta
$$
uniformly for $1\leq h\leq X$.
}

\medskip

\noindent{\bf Remark 2.}
For holomorphic cusp forms, Theorem 2 improves the result in Luo \cite{L}, where it is proved that
$$
\sum_{n\geq 1}\lambda_f(n+h)r(n)\phi\left(\frac{n}{X}\right)
\ll_{f,\varepsilon}X^{1-\frac{1}{4}+\varepsilon}\Delta^2.
$$

Assume that $\lambda_f(n)\ll n^{\theta+\varepsilon}$. By Theorem 2 we have
$$
\sum_{X/2\leq n\leq X}\lambda_f(n+h)r(n)=\sum_{n\geq 1}\lambda_f(n+h)r(n)\phi\left(\frac{n}{X}\right)+
O_{\varepsilon}\left(X^{1+\theta+\varepsilon}\Delta^{-1}\right)
\ll_{f,\varepsilon}X^{\frac{3}{4}+\varepsilon}\Delta+X^{1+\theta+\varepsilon}\Delta^{-1}.
$$
On taking $\Delta=X^{\frac{1}{8}+\frac{\theta}{2}}$, we obtain
\bna
\sum_{X/2\leq n\leq X}\lambda_f(n+h)r(n)\ll_{f,\varepsilon} X^{1-\frac{1-4\theta}{8}+\varepsilon},
\ena
uniformly for $1\leq h\leq X$. Thus we get the following result which improves (1.2).

\medskip

\noindent {\bf Corollary 1.}\, {\it Let $f$ be a holomorphic cusp form of weight $\kappa$ or
a Maass cusp form of
Laplace eigenvalue $\frac{1}{4}+\mu^2$ for $SL(2,\mathbb{Z})$ and
denote by $\lambda_f(n)$ its $n$-th Hecke eigenvalue.
Assume that $\lambda_f(n)\ll n^{\theta+\varepsilon}$.
Then for $X>1$ and any $\varepsilon>0$, we have
\bna
\sum_{n\leq X}\lambda_f(n+h)r(n)\ll_{f,\varepsilon}  X^{1-\frac{1-4\theta}{8}+\varepsilon}
\ena
uniformly for $1\leq h\leq X$.
}

\medskip

For the proof of Theorem 2, we use a different method from \cite{L}, \cite{LWZ} and Theorem 1.
Precisely, as in \cite{S}, we apply the classical Hardy-Littlewood-Kloosterman circle method to transform the
sum in Theorem 2. Then we can apply the Voronoi formula for the cusp form $f$ and
an asymptotic formula for the sum $\sum_{|m|\leq \sqrt{X}}e(\alpha m^2)$, $\alpha\in [0,1]$,
to get essentially three sums involving transforms of Bessel functions and some
exponential sums. Finally, Theorem 2 follows from nontrivial estimates of these
transforms and exponential sums.

\medskip

\section{Voronoi formulas}
\setcounter{equation}{0}
\medskip

We need Voronoi formulas for both $\lambda_f(n)$
and $r_{\ell}(n)$. Let $J_s$, $K_s$ and $Y_s$ denote the standard $J$-Bessel function,
$K$-Bessel function and $Y$-Bessel function, respectively.
Let $f$ be a general $GL(2)$ cusp form, that is (as in \cite{BHM}) a holomorphic cusp
form of integral weight $\kappa\geq 1$, level $D$ and nebentypus $\chi_{D}$, or a
Maass cusp form of weight $0$ or $1$, level $D$, nebentypus $\chi_D$ and Laplace eigenvalue
$1/4+\mu^2$. Suppose that $f$ is a newform and denote by $\lambda_f(n)$
its $n$-th Hecke eigenvalue.
Then we have the following  Voronoi formula for $\lambda_f(n)$ (see \cite{BHM} or \cite{HM}, Proposition 2.1)

\medskip

\noindent{\bf Lemma 1.} {\it Let $f$ be a cusp form (holomorphic or Maass)
of weight $\kappa$, level $D$ and nebentypus $\chi_D$. Let $q\equiv 0 \bmod D$ and
let $a$ be an integer coprime to $q$ and $ad\equiv1(\bmod q)$.
If $v(x)$ is a compactly supported smooth function on $\mathbb{R}^+$, then
\bna
\sum_{n\geq 1}\lambda_f(n)e\left(\frac{an}{q}\right)v(n)
&=&\frac{\chi_D(d)}{q}\sum_{n\geq 1}\lambda_f(n)
e\left(-\frac{dn}{q}\right)\mathcal {V}^{+}\left(\frac{n}{q^2}\right)\\
&&+\omega_f
\frac{\Gamma\left(\frac{1}{2}+i\mu-\frac{\kappa}{2}\right)}{\Gamma\left(\frac{1}{2}+i\mu+\frac{\kappa}{2}\right)}
\frac{\chi_D(d)}{q}\sum_{n\geq 1}\lambda_f(n)
e\left(\frac{dn}{q}\right)\mathcal {V}^{-}\left(\frac{n}{q^2}\right),
\ena
where $\omega_f$ is a constant depending only on $f$, $\mu$ is the spectral parameter of $f$
in the Maass case, and
$$
\mathcal {V}^{\pm}(y)=\int_{0}^{\infty}v(x)\mathcal {H}_f^{\pm}(4\pi\sqrt{xy})dx.
$$
In this formula, (i) if $f$ is induced from a holomorphic form of weight $k$,
\bea
\mathcal {H}_f^+(x)=2\pi i^k J_{k-1}(x),\qquad \mathcal {H}_f^-(x)=0;
\eea
(ii) if $\kappa$ is even, $f$ is not induced from a holomorphic form,
\bea
\mathcal {H}_f^+(x)=\frac{-\pi}{\cosh(\pi \mu)}\left\{Y_{2i\mu}(x)+Y_{-2i\mu}(x)\right\},\qquad
\mathcal {H}_f^-(x)=4\cosh(\pi \mu)K_{2i\mu}(x);
\eea
(iii) if $\kappa$ is odd, $f$ is not induced from a holomorphic form,
\bea
\mathcal {H}_f^+(x)=\frac{\pi}{\sinh(\pi \mu)}\left\{Y_{2i\mu}(x)-Y_{-2i\mu}(x)\right\},\qquad
\mathcal {H}_f^-(x)=-4i\sinh(\pi \mu)K_{2i\mu}(x).
\eea
}

$\mathcal {V}^{\pm}(y)$ has the following properties.

\noindent{\bf Lemma 2.} {\it If $v(x)$ is a smooth function supported on $[AY,BY]$ ($B>A>0$)
satisfying $x^jv^{(j)}(x)\ll_{A,B,j}1 $, then for any integer $j\geq 0$ and any $\varepsilon>0$,
$$
\mathcal {V}^{\pm}(y)\ll_{f,A,B,j,\varepsilon}Y\left(1+\sqrt{yY}\right)^{-\frac{1}{2}}
\left(1+\frac{1}{\sqrt{yY}}\right)^{2|\Im \mu|+\varepsilon}
\left(\frac{1}{\sqrt{yY}}+\frac{1}{yY}\right)^{j}
$$
where $\mu$ is the spectral parameter of $f$
in the Maass form case and $\mu=0$ if $f$ is a holomorphic form.
}

\medskip

\noindent {\bf Proof.} We follow closely \cite{BHM} (see Section 4.3 and Appendix 1).
Set $R=4\pi \sqrt{y}$. Let $H_s$ denote either $J_s$, $Y_{s}$ or $K_s$. Then $H_s$ satisfies the
recurrence relations (see (6.1) in \cite{BHM}),
\bea
\frac{\mathrm{d}}{\mathrm{d}x}\left\{(R\sqrt{x})^{s+1}H_{s+1}(R\sqrt{x})\right\}
=\pm \frac{R^2}{2}(R\sqrt{x})^sH_s\left(R\sqrt{x}\right)，
\eea
where the $-$ appears only for $H_s=K_s$. By the definition of $\mathcal {V}^{\pm}(y)$,
we only need to prove that
\bea
\mathcal {V}(y)&:=&\int_0^{\infty}v(x)H_s\left(R\sqrt{x}\right)\mathrm{d}x\nonumber\\
&&\ll_{f,A,B,j,\varepsilon}
Y\left(1+R\sqrt{Y}\right)^{-\frac{1}{2}}\left(1+\frac{1}{R\sqrt{Y}}\right)^{2|\Im \mu|+\varepsilon}
\left(\frac{1}{R\sqrt{Y}}+\frac{1}{R^2Y}\right)^{j}.
\eea
By integration by parts $j$ times with respect to $x$ and (2.4), we have
\bea
\mathcal {V}(y)=\left(\mp \frac{R^2}{2}\right)^{-j}\int_0^{\infty}
\frac{\mathrm{d}^j}{\mathrm{d}x^j}\left\{v(x)\left(R\sqrt{x}\right)^{-s}\right\}
\left(R\sqrt{x}\right)^{s+j}H_{s+j}\left(R\sqrt{x}\right)\mathrm{d}x.
\eea
Note that
\bna
\frac{\mathrm{d}^j}{\mathrm{d}x^j}\left\{v(x)\left(R\sqrt{x}\right)^{-s}\right\}
=R^{-s}\sum_{i=0}^jC_{j}^iv^{(i)}(x)\left(-\frac{s}{2}\right)
\left(-\frac{s}{2}-1\right)\cdot\cdot\cdot \left(-\frac{s}{2}-j+i+1\right)
x^{-\frac{s}{2}-(j-i)}
\ena
and by Proposition 6.2 in \cite{BHM}, for any $\varepsilon>0$,
\bna
H_{s+j}\left(R\sqrt{x}\right)\ll_{f,j,\varepsilon}\left(1+R\sqrt{x}\right)^{-\frac{1}{2}}
\left(1+\frac{1}{R\sqrt{x}}\right)^{j+2|\Im \mu|+\varepsilon}.
\ena
Plugging these estimates into (2.6) we obtain (2.5). This proves Lemma 2. \hfill $\Box$

\medskip

The following Voronoi formula for $r_{\ell}(n)$ is due to Luo (see \cite{L}, Lemma 2).

\noindent{\bf Lemma 3.} {\it Let $a$ be an integer coprime to $q$, $ad\equiv1(\bmod q)$ and
\bea
\epsilon_d=\left\{\begin{array}{ll}1,&\mbox{if\, $d\equiv 1(\bmod 4)$},\\
i,&\mbox{if \,$d\equiv -1(\bmod 4)$.}
\end{array}\right.\eea
If $w(x)$ is a compactly supported smooth function on $\mathbb{R}^+$, then
\bna
\sum_{n\geq 1}r_{\ell}(n)e\left(\frac{an}{q}\right)w(n)
&=&\left(\frac{2\pi i}{q}\right)^{\frac{\ell}{2}}
\Gamma\left(\frac{\ell}{2}\right)^{-1}
\left(\left(\frac{q}{d}\right)\epsilon_d^{-1}\right)^{\ell}
\tilde{w}\left(\frac{\ell}{2}\right)\nonumber\\
&&+\frac{2\pi i^{\frac{\ell}{2}}}{q}\left(\left(\frac{q}{d}\right)\epsilon_d^{-1}\right)^{\ell}
\sum_{n\geq 1}r_{\ell}(n)
e\left(-\frac{dn}{q}\right)n^{\frac{1-\frac{\ell}{2}}{2}}\mathcal {W}\left(\frac{n}{q^2}\right),\nonumber\\
\ena
where $\tilde{w}(s)=\int_0^{\infty}w(x)x^{s-1}\mathrm{d}x$ and
$$
\mathcal {W}(y)=\int_{0}^{\infty}w(x)x^{\frac{\frac{\ell}{2}-1}{2}}J_{\frac{\ell}{2}-1}(4\pi\sqrt{xy})dx.
$$
}

$\mathcal {W}(y)$ has the following properties.

\noindent{\bf Lemma 4.} {\it If $w(x)$ is a smooth function supported on $[AY,BY]$ ($B>A>0$)
satisfying $w^{(j)}(x)\ll_{A,B,j} \left(\Delta/Y\right)^j$ ($\log \Delta \asymp \log Y$),
then for any $s\geq 0$ and any integer $j\geq 0$,
$$
\mathcal {W}_s(y)=\int_{0}^{\infty}w(x)x^{\frac{s}{2}}J_s(4\pi\sqrt{xy})dx
\ll_{s,A,B,j}Y^{1+\frac{s}{2}}\left(1+\sqrt{yY}\right)^{-\frac{1}{2}}\left(1+\frac{1}{\sqrt{yY}}\right)^{-s}
\left(\frac{\Delta}{1+\sqrt{yY}}\right)^{j}.
$$
}

\medskip

\noindent{\bf Proof.} The proof is very similar as that of Lemma 2 and one uses the estimate
(see Proposition 6.1 in \cite{BHM})
\bna
J_{s}\left(R\sqrt{x}\right)\ll_{s}\left(1+R\sqrt{x}\right)^{-\frac{1}{2}}
\left(1+\frac{1}{R\sqrt{x}}\right)^{-s}
\ena
for any $s\geq 0$.

\medskip

For the properties of Bessel functions, we quote the following lemma
(see \cite{GR}, page 920, 8.451-1, 8.451-2 and 8.451-6).

\medskip

\noindent {\bf Lemma 5.} {\it For $s\in \mathbb{C}$ and $x\gg 1$, we have
\bna
&&J_{s}(x)=\sqrt{\frac{2}{\pi x}}\left\{\cos\left(x-\frac{\pi}{2}s-\frac{\pi}{4}\right)
-\frac{4{s}^2-1}{8x}\sin\left(x-\frac{\pi}{2}s-\frac{\pi}{4}\right)
+O_{s}\left(x^{-2}\right)\right\},\\
&&Y_{s}(x)=\sqrt{\frac{2}{\pi x}}\left\{\sin\left(x-\frac{\pi s}{2}-\frac{\pi}{4}\right)+
\frac{4s^2-1}{8x}\cos\left(x-\frac{\pi s}{2}-\frac{\pi}{4}\right)
+O_s\left(x^{-2}\right)\right\},\\
&&K_{s}(x)=\sqrt{\frac{\pi}{2 x}}e^{-x}\left\{1+O_{s}\left(x^{-1}\right)\right\}.
\ena
}

\section{Proof of Theorem 1}
\setcounter{equation}{0}
\medskip

Let $\phi(x)$ be a smooth function which is supported on $[1/2,1]$ and identically equal to 1 on
$[1/2+\Delta^{-1},1-\Delta^{-1}]$ with $\log \Delta\asymp \log X$ to be chosen later, satisfying
$\phi^{(j)}(x)\ll_j \Delta^{j}$ for all integers $j\geq 0$. Then
\begin{equation}
\sum_{X/2<n\leq X}\lambda_f(n+h)r_{\ell}(n)=
\sum_{n\geq 1}\lambda_f(n+h)r_{\ell}(n)\phi\left(\frac{n}{X}\right)+E_h(X),
\end{equation}
where
\bea
E_h(X)=\sum_{X/2<n\leq X}\lambda_f(n+h)r_{\ell}(n)\left(1-\phi\left(\frac{n}{X}\right)\right)
\ll_{\ell,\varepsilon} X^{\frac{\ell}{2}+\theta+\varepsilon}\Delta^{-1},
\eea
using the bound $r_{\ell}(n)\ll_{\ell,\varepsilon} n^{\frac{\ell}{2}-1+\varepsilon}$ and
$\lambda_f(n)\ll n^{\theta+\varepsilon}$ for any $\varepsilon>0$.

Denote the first term on the right side of (3.1) by $\mathcal{S}_h(X)$. Then we have
\bna
\mathcal{S}_h(X)=\sum_{n\geq 1}\lambda_f(n+h)r_{\ell}(n)\phi\left(\frac{n}{X}\right)
\varphi\left(\frac{n+h}{X+h}\right),
\ena
where $\varphi(x)$ is a smooth function supported on $[1/4,5/4]$ and equals 1 if $x\in [1/2,1]$,
satisfying $\varphi^{(j)}(x)\ll_j 1$. Let $\delta(i,j)=\left\{\begin{array}{ll}
1,&\mbox{if $i=j$,}\\0,&\mbox{otherwise}.\end{array}
\right.$
Then
\bna
\delta(i,j)=\int_{0}^1e\left((i-j)\beta\right)\mathrm{d}\beta.
\ena
and
\bna
\mathcal{S}_h(X)&=&\sum_{n\geq 1}r_{\ell}(n)\phi\left(\frac{n}{X}\right)\sum_{m\geq 1}\lambda_f(m)
\varphi\left(\frac{m}{X+h}\right)\delta(n+h,m)\\
&=&\int_0^1e(h\beta)\sum_{n\geq 1}r_{\ell}(n)\phi\left(\frac{n}{X}\right)e(\beta n)
\sum_{m\geq 1}\lambda_f(m)
\varphi\left(\frac{m}{X+h}\right)e(-\beta m)\mathrm{d}\beta.
\ena
Applying Jutila's circle method (see \cite{J1}), we approximate $\delta(i,j)$ by
\bna
\widetilde{\delta}(i,j)=\int_{0}^1\widetilde{I}_{\mathcal{Q},\delta}(\beta)
e\left((i-j)\beta\right)\mathrm{d}\beta,
\ena
and $\mathcal{S}_h(X)$ by
\bea
\widetilde{\mathcal{S}}_h(X)
=\int_0^1\widetilde{I}_{\mathcal{Q},\delta}(\beta)e(h\beta)
\sum_{n\geq 1}r_{\ell}(n)\phi\left(\frac{n}{X}\right)e(\beta n)\sum_{m\geq 1}\lambda_f(m)
\varphi\left(\frac{m}{X+h}\right)e(-\beta m)\mathrm{d}\beta,
\eea
where $\widetilde{I}_{\mathcal{Q},\delta}(\beta)$ is defined in (1.3),
$Q>0$ is a parameter to be chosen soon, $\mathcal{Q}\subset [1,Q]$, $L=\sum_{q\in \mathcal{Q}}\phi(q)$
and $Q^{-2}\leq \delta \leq Q^{-1}$.

By Cauchy's inequality and (1.4), we have, for $L\gg_{\varepsilon} Q^{2-\varepsilon}$,
\bna
&&|\mathcal{S}_h(X)-\widetilde{\mathcal{S}}_h(X)|\\
&\leq& \int_0^1\left|1-\widetilde{I}_{\mathcal{Q},\delta}(\beta)\right|
\left|\sum_{n\geq 1}r_{\ell}(n)\phi\left(\frac{n}{X}\right)e(\beta n)\right|
\left|\sum_{m\geq 1}\lambda_f(m)
\varphi\left(\frac{m}{X+h}\right)e(-\beta m)\right|\mathrm{d}\beta\\
&\ll_{f,\varepsilon}& X^{\frac{1}{2}+\varepsilon}
\left(\int_0^1\left|1-\widetilde{I}_{\mathcal{Q},\delta}(\beta)\right|^2\mathrm{d}\beta\right)^{1/2}
\left(\int_0^1\left|\sum_{n\geq 1}r_{\ell}(n)\phi\left(\frac{n}{X}\right)e(\beta n)\right|^2
\mathrm{d}\beta\right)^{1/2}\\
&\ll_{f,\ell,\varepsilon}&\frac{X^{\frac{\ell}{2}+\varepsilon}}{\sqrt{\delta}Q},
\ena
where we have used the trivial bound
$r_{\ell}(n)\ll_{\ell,\varepsilon} n^{\frac{\ell}{2}-1+\varepsilon}$ and
the estimate (see (4.2) in \cite{BHM})
\[
\sum_{n\leq x}\lambda_f(n)e(\alpha n)\ll_f x^{\frac{1}{2}}\log 2x
\]
which is uniform in $\alpha\in \mathbb{R}$. Take $\delta=X^{-1}$. We have
\bea
\mathcal{S}_h(X)=\widetilde{\mathcal{S}}_h(X)+O_{f,\ell,\varepsilon}\left(X^{\frac{\ell+1}{2}+\varepsilon}Q^{-1}\right).
\eea

In the following, we estimate $\widetilde{\mathcal{S}}_h(X)$ in (3.3). By (1.3), we have
\bea
\widetilde{\mathcal{S}}_h(X)=\frac{1}{2\delta}\int_{-\delta}^{\delta}
\widetilde{\mathcal{S}}_h(X,\beta)e(h \beta)\mathrm{d}\beta,
\eea
where
\bea
\widetilde{\mathcal{S}}_h(X,\beta)&=&
\frac{1}{L}\sum_{q\in \mathcal{Q}}\sum_{a=1 \atop (a,q)=1}^qe\left(\frac{ah}{q}\right)
\sum_{n\geq 1}r_{\ell}(n)e\left(\frac{an}{q}\right)\phi\left(\frac{n}{X}\right)e(\beta n)\nonumber\\
&&\qquad\qquad \qquad \qquad
\times\sum_{m\geq 1}\lambda_f(m)e\left(-\frac{am}{q}\right)\varphi\left(\frac{m}{X+h}\right)e(-\beta m).
\eea

Now we choose the set of moduli $\mathcal{Q}$ as follows
\bea
\mathcal{Q}=\{4Dp: p\in [Q/(8D),Q/(4D)] \ \text{is prime and}\ (p,2Dh)=1\}.
\eea
Then the requirement $L\gg_{\varepsilon} Q^{2-\varepsilon}$ is satisfied.
Applying Lemma 1 with $v(x)=\varphi\left(\frac{x}{X+h}\right)e(-\beta x)$ to
the $m$-sum in (3.6), we have
\bna
&&\sum_{m\geq 1}\lambda_f(m)e\left(-\frac{am}{q}\right)\varphi\left(\frac{m}{X+b}\right)e(-\beta m)\\
&=&\frac{\chi_D(-d)}{q}\sum_{m\geq 1}\lambda_f(m)
e\left(\frac{dm}{q}\right)\mathcal {V}_{\beta}^{+}\left(\frac{m}{q^2}\right)\\
&&+\omega_f
\frac{\Gamma\left(\frac{1}{2}+i\mu-\frac{\kappa}{2}\right)}{\Gamma\left(\frac{1}{2}
+i\mu+\frac{\kappa}{2}\right)}
\frac{\chi_D(-d)}{q}\sum_{m\geq 1}\lambda_f(m)
e\left(-\frac{dm}{q}\right)\mathcal {V}_{\beta}^{-}\left(\frac{m}{q^2}\right),
\ena
where $ad\equiv 1(\bmod q)$ and
\bna
\mathcal {V}_{\beta}^{\pm}(y)=\int_0^{\infty}\varphi\left(\frac{x}{X+h}\right)
e(-\beta x)\mathcal {H}_f^{\pm}(4\pi\sqrt{xy})\mathrm{d}x
\ena
with $\mathcal {H}_f^{\pm}(x)$ defined in (2.1)-(2.3). Note that for $|\beta|\leq \delta=X^{-1}$,
$x^jv^{(j)}(x)\ll_j 1$. By Lemma 2, the contribution from
$mX/q^2\gg X^{\varepsilon}$ is negligible. Thus
\bea
\widetilde{\mathcal{S}}_h(X,\beta)=\widetilde{\mathcal{S}}_h(X,\beta,+)
+\omega_f
\frac{\Gamma\left(\frac{1}{2}+i\mu-\frac{\kappa}{2}\right)}{\Gamma\left(\frac{1}{2}
+i\mu+\frac{\kappa}{2}\right)}\widetilde{\mathcal{S}}_h(X,\beta,-)+O_{\varepsilon}(1),
\eea
where
\bea
\widetilde{\mathcal{S}}_h(X,\beta,\pm)&=&
\frac{1}{L}\sum_{q\in \mathcal{Q}}\frac{1}{q}\sum_{a=1 \atop (a,q)=1}^q
\chi_D(-d)e\left(\frac{ah}{q}\right)
\sum_{1\leq m\ll Q^2X^{\varepsilon}/X}\lambda_f(m)
e\left(\pm\frac{dm}{q}\right)\mathcal {V}_{\beta}^{\pm}\left(\frac{m}{q^2}\right)\nonumber\\
&&\sum_{n\geq 1}r_{\ell}(n)e\left(\frac{an}{q}\right)\phi\left(\frac{n}{X}\right)e(\beta n).
\eea

Applying Lemma 3 with $w(x)=\phi\left(\frac{x}{X}\right)e(\beta x)$ to the $n$-sum in (3.9), we have
\bea
\sum_{n\geq 1}r_{\ell}(n)e\left(\frac{an}{q}\right)\phi\left(\frac{n}{X}\right)e(\beta n)
&=&\left(\frac{2\pi i}{q}\right)^{\frac{\ell}{2}}
\Gamma\left(\frac{\ell}{2}\right)^{-1}
\left(\left(\frac{q}{d}\right)\epsilon_d^{-1}\right)^{\ell}\tilde{w}_{\beta}\left(\frac{\ell}{2}\right)\nonumber\\
&&+\frac{2\pi i^{\frac{\ell}{2}}}{q}\left(\left(\frac{q}{d}\right)\epsilon_d^{-1}\right)^{\ell}
\sum_{n\geq 1}r_{\ell}(n)
e\left(-\frac{dn}{q}\right)n^{\frac{1-\frac{\ell}{2}}{2}}\mathcal {W}_{\beta}\left(\frac{n}{q^2}\right),\nonumber\\
\eea
where
\bea
 \tilde{w}_{\beta}(s)=\int_0^{\infty}\phi\left(\frac{x}{X}\right)e(\beta x)x^{s-1}\mathrm{d}x,
\eea
and
\bea
\mathcal {W}_{\beta}(y)=\int_{0}^{\infty}\phi\left(\frac{x}{X}\right)e(\beta x)x^{\frac{\frac{\ell}{2}-1}{2}}J_{\frac{\ell}{2}-1}(4\pi\sqrt{xy})dx.
\eea
Note that for $|\beta|\leq \delta=X^{-1}$,
$w^{(j)}(x)\ll_j (\Delta/X)^j$.
By Lemma 4, the contribution from $nX/q^2\gg \Delta^2X^{\varepsilon}$
in the second sum in (3.10) is negligible.

Plugging (3.10) into (3.9) we obtain
\bea
\widetilde{\mathcal{S}}_h(X,\beta,\pm)=\widetilde{\mathcal{S}}_h^0(X,\beta,\pm)
+\widetilde{\mathcal{S}}_h^{\flat}(X,\beta,\pm)+O_{\varepsilon}(1),
\eea
where
\bna
\widetilde{\mathcal{S}}_h^0(X,\beta,\pm)&=&
\frac{(2\pi i)^{\frac{\ell}{2}}}{\Gamma\left(\ell/2\right)}
\frac{\tilde{w}_{\beta}\left(\ell/2\right)}{L}
\sum_{q\in \mathcal{Q}}\frac{1}{q^{\frac{\ell}{2}+1}}
\sum_{1\leq m\ll Q^2X^{\varepsilon}/X}\lambda_f(m)
\mathcal {V}_{\beta}^{\pm}\left(\frac{m}{q^2}\right)\\
&&\qquad\qquad\times
\sum_{a=1 \atop ad\equiv 1(\bmod q)}^q\chi_D(-d)\left(\left(\frac{q}{d}\right)\epsilon_d^{-1}\right)^{\ell}
e\left(\frac{ha\pm md}{q}\right),\\
\widetilde{\mathcal{S}}_h^{\flat}(X,\beta,\pm)&=&
\frac{2\pi i^{\frac{\ell}{2}}}{L}\sum_{q\in \mathcal{Q}}\frac{1}{q^2}
\sum_{1\leq m\ll Q^2X^{\varepsilon}/X}\lambda_f(m)
\mathcal {V}_{\beta}^{\pm}\left(\frac{m}{q^2}\right)
\sum_{1\leq n \ll \Delta^2Q^2X^{\varepsilon}/X}r_{\ell}(n)
n^{\frac{1-\frac{\ell}{2}}{2}}\mathcal {W}_{\beta}\left(\frac{n}{q^2}\right)\\
&&\times
\sum_{a=1 \atop ad\equiv 1(\bmod q)}^q\chi_D(-d)\left(\left(\frac{q}{d}\right)\epsilon_d^{-1}\right)^{\ell}
e\left(\frac{ha+(\pm m-n)d}{q}\right).
\ena

\medskip

{\bf Bounding the character sum.} (i) $\ell$ is even. In this case
$\left(\left(\frac{q}{d}\right)\epsilon_d^{-1}\right)^{\ell}
=\left(\left(\frac{-q}{a}\right)\epsilon_a\right)^{\ell}=\chi_4(a)^{\frac{\ell}{2}}$.
By our choice of $\mathcal{Q}$ in (3.7), we have, for $M\in \mathbb{Z}$,
\bea
&&\sum_{a=1 \atop ad\equiv 1(\bmod q)}^q\chi_D(-d)\left(\left(\frac{q}{d}\right)
\epsilon_d^{-1}\right)^{\ell}e\left(\frac{ha+Md}{q}\right)\nonumber\\
&=&\sum_{a_1\bmod 4D \atop a_1\bar{a}_1\equiv 1(\bmod 4D)}
\chi_4(a_1p)^{\frac{\ell}{2}}\chi_{D}(-\bar{a}_1\overline{p})
e\left(\frac{ha_1+\bar{p}^2M\bar{a}_1}{4D}\right)
\sum_{a_2\bmod p \atop a_2\bar{a}_2\equiv 1(\bmod p)}
e\left(\frac{ha_2+\overline{4D}^2M\bar{a}_2}{p}\right)\nonumber\\
&\ll_D& \sqrt{p}
\eea
by Weil's bound for Kloosterman sum.

(ii) $\ell$ is odd. In this case
$\left(\left(\frac{q}{d}\right)\epsilon_d^{-1}\right)^{\ell}
=\left(\left(\frac{-q}{a}\right)\epsilon_a\right)^{\ell}
=\left(\frac{-q}{a}\right)\epsilon_a^{\ell}=\left(\frac{q}{a}\right)\chi_4(a)\epsilon_a^{\ell}$.
By our choice of $\mathcal{Q}$ in (3.7) and quadratic reciprocity, we have, for $M\in \mathbb{Z}$,
\bea
&&\sum_{a=1 \atop ad\equiv 1(\bmod q)}^q\chi_D(-d)\left(\left(\frac{q}{d}\right)
\epsilon_d^{-1}\right)^{\ell}e\left(\frac{ha+Md}{q}\right)\nonumber\\
&=&\sum_{a_1\bmod 4D \atop a_1\bar{a}_1\equiv 1(\bmod 4D)}\left(\frac{pa_1}{4D}\right)
(-1)^{\frac{p-1}{2}\frac{pa_1-1}{2}}\chi_4(pa_1)\chi_{D}(-\bar{a}_1\overline{p})
\epsilon_{pa_1}^{\ell}e\left(\frac{ha_1+\bar{p}^2M\bar{a}_1}{4D}\right)\nonumber\\
&&\sum_{a_2\bmod p \atop a_2\bar{a}_2\equiv 1(\bmod p)}\left(\frac{4Da_2}{p}\right)
e\left(\frac{ha_2+\overline{4D}^2M\bar{a}_2}{p}\right)\nonumber\\
&\ll_D& \sqrt{p}
\eea
by the well known bound for Sali\'{e} sums (see Corollary 4.10 in \cite{I2}).

\medskip

{\bf Bounding $\widetilde{\mathcal{S}}_h^0(X,\beta,\pm)$.} Denote $\vartheta=|\Im \mu|$. Then
by Selberg's bound $\vartheta\leq 1/4$. By Lemma 2, for $p\in [Q/(8D),Q/(4D)]$, we have
\bea
\sum_{1\leq m\ll Q^2X^{\varepsilon}/X}\lambda_f(m)
\mathcal {V}_{\beta}^{\pm}\left(\frac{m}{q^2}\right)&\ll&
X\sum_{1\leq m\ll Q^2X^{\varepsilon}/X}\left|\lambda_f(m)\right|
\left(1+\sqrt{\frac{16D^2p^2}{mX}}\right)^{2\vartheta+\varepsilon}\nonumber\\
&\ll_f&X\sum_{1\leq \frac{mX}{16D^2p^2}\ll X^{\varepsilon}}\left|\lambda_f(m)\right|
+X\sum_{\frac{mX}{16D^2p^2}\leq 1}\left|\lambda_f(m)\right|
\left(\frac{Q^2}{mX}\right)^{\vartheta+\varepsilon}\nonumber\\
&\ll_f&Q^2X^{\varepsilon}
\eea
by the Rankin-Selberg estimate
\bea
\sum_{m\leq x}\left|\lambda_f(m)\right|^2\ll_f x.
\eea

Recall that $L\gg_{\varepsilon}Q^{2-\varepsilon}$.
Bounding the integral $\tilde{w}_{\beta}\left(\frac{\ell}{2}\right)$ in (3.11)
trivially and by (3.14)-(3.16), we have
\bea
\widetilde{\mathcal{S}}_h^0(X,\beta,\pm)&\ll_{f,\ell,\varepsilon}&
\frac{X^{\frac{\ell}{2}}}{L}\sum_{p\in [Q/(8D),Q/(4D)]}\frac{1}{p^{(\ell+1)/2}}
\sum_{1\leq m\ll Q^2X^{\varepsilon}/X}\left|\lambda_f(m)\right|
\left|\mathcal {V}^{\pm}\left(\frac{m}{16D^2p^2}\right)\right|\nonumber\\
&\ll_{f,\ell,\varepsilon}&\frac{X^{\frac{\ell}{2}+\varepsilon}}{Q^{\frac{\ell-1}{2}}}.
\eea

{\bf Bounding $\widetilde{\mathcal{S}}_h^{\flat}(X,\beta,\pm)$.}
First, for $n \ll \Delta^2Q^2X^{\varepsilon}/X$,
we estimate $\mathcal {W}_{\beta}\left(n/q^2\right)$ more precisely. By (3.12) we have
\bna
\mathcal {W}_{\beta}\left(\frac{n}{q^2}\right)=X^{\frac{\ell}{4}+\frac{1}{2}}
\int_{0}^{\infty}\phi\left(x\right)e(\beta Xx)x^{\frac{\frac{\ell}{2}-1}{2}}
J_{\frac{\ell}{2}-1}\left(\frac{4\pi\sqrt{nXx}}{q}\right)dx.
\ena
Set $R=\frac{4\pi\sqrt{nX}}{q}$. By partial integration once with respect to $x$ and
applying the recurrence relation $(x^sJ_s(x))'=x^sJ_{s-1}(x)$,
we have
\bea
\mathcal {W}_{\beta}\left(\frac{n}{q^2}\right)=X^{\frac{\ell}{4}+\frac{1}{2}}R^{-\frac{\ell}{2}+1}\left(\frac{-2}{R^2}\right)
\int_{0}^{\infty}\frac{d}{dx}\left\{\phi\left(x\right)e(\beta Xx)\right\}
\left(R\sqrt{x}\right)^{\frac{\ell}{2}}
J_{\frac{\ell}{2}}\left(R\sqrt{x}\right)dx:=\mathcal{E}_1+\mathcal{E}_2,
\eea
where
\bna
\mathcal{E}_1&=&-4\pi i \beta X^{\frac{\ell}{4}+\frac{3}{2}}R^{-\frac{\ell}{2}-1}
\int_{0}^{\infty}\phi(x) e(\beta Xx)\left(R\sqrt{x}\right)^{\frac{\ell}{2}}
J_{\frac{\ell}{2}}\left(R\sqrt{x}\right)dx,\\
\mathcal{E}_2&=&-2X^{\frac{\ell}{4}+\frac{1}{2}}R^{-\frac{\ell}{2}-1}
\int_{0}^{\infty}\phi'(x)e(\beta Xx)\left(R\sqrt{x}\right)^{\frac{\ell}{2}}
J_{\frac{\ell}{2}}\left(R\sqrt{x}\right)dx.
\ena
Using the bound $J_s(x)\ll \min\{x^s, x^{-\frac{1}{2}}\}$ for $s\geq 0$, we have, for $|\beta|\leq X^{-1}$,
\bea
\mathcal{E}_1\ll X^{\frac{\ell}{4}+\frac{1}{2}}R^{-1}\min\{R^{\frac{\ell}{2}}, R^{-1/2}\},
\eea
and
\bea
\mathcal{E}_2&=&-2X^{\frac{\ell}{4}+\frac{1}{2}}R^{-\frac{\ell}{2}-1}
\left\{\int_1^{1+2\Delta^{-1}}+
\int_{2-2\Delta^{-1}}^2\right\}\phi'(x)e(\beta Xx)\left(R\sqrt{x}\right)^{\frac{\ell}{2}}
J_{\frac{\ell}{2}}\left(R\sqrt{x}\right)dx\nonumber\\
&\ll&X^{\frac{\ell}{4}+\frac{1}{2}}R^{-1}\min\{R^{\frac{\ell}{2}}, R^{-1/2}\}.
\eea
By (3.19)-(3.21), we have
\bea
\mathcal {W}_{\beta}\left(\frac{n}{q^2}\right)
\ll X^{\frac{\ell}{4}+\frac{1}{2}}R^{-1}\min\{R^{\frac{\ell}{2}}, R^{-1/2}\}
\ll X^{\frac{\ell}{4}+\frac{1}{2}}
\min\left\{\left(\frac{\sqrt{nX}}{q}\right)^{\frac{\ell}{2}-1},
\left(\frac{q}{\sqrt{nX}}\right)^{\frac{3}{2}}\right\}.
\eea

By (3.16), (3.22) and the estimate $r_{\ell}(n)\ll n^{\frac{\ell}{2}-1+\varepsilon}$, we have
\bea
\widetilde{\mathcal{S}}_h^{\flat}(X,\beta,\pm)&\ll&
\sum_{p\in [Q/(8D),Q/(4D)]}\frac{1}{p^{3/2}}
\sum_{n\ll \Delta^2Q^2X^{\varepsilon}/X}n^{\frac{\ell}{4}-\frac{1}{2}}
X^{\frac{\ell}{4}+\frac{1}{2}}
\min\left\{\left(\frac{\sqrt{nX}}{p}\right)^{\frac{\ell}{2}-1},
\left(\frac{p}{\sqrt{nX}}\right)^{\frac{3}{2}}\right\}\nonumber\\
&\ll&\frac{X^{\frac{\ell}{4}+\frac{1}{2}+\varepsilon}}{\sqrt{Q}}
\sum_{n\ll \Delta^2Q^2X^{\varepsilon}/X}n^{\frac{\ell}{4}-\frac{1}{2}}
\min\left\{\left(\frac{\sqrt{nX}}{Q}\right)^{\frac{\ell}{2}-1},
\left(\frac{Q}{\sqrt{nX}}\right)^{\frac{3}{2}}\right\}\nonumber\\
&\ll&\frac{X^{\frac{\ell}{4}+\frac{1}{2}+\varepsilon}}{\sqrt{Q}}\left\{
\sum_{Q^2/X<n\ll Q^2\Delta^2X^{\varepsilon}/X}\left(\frac{Q}{\sqrt{nX}}\right)^{\frac{3}{2}}
+\sum_{n\leq Q^2/X}\left(\frac{\sqrt{nX}}{Q}\right)^{\frac{\ell}{2}-1}\right\}\nonumber\\
&\ll&X^{\frac{\ell}{4}+\frac{1}{2}+\varepsilon}Q^{-\frac{1}{2}}
+X^{\frac{\ell}{4}-\frac{1}{2}+\varepsilon}Q^{\frac{3}{2}}\Delta^{\frac{1}{2}}.
\eea

By (3.4), (3.5), (3.8), (3.13), (3.18) and (3.23), we have, for $Q<X$,
\bea
\mathcal{S}_h(X)\ll\frac{X^{\frac{\ell}{2}+\varepsilon}}{Q^{\frac{\ell-1}{2}}}
+X^{\frac{\ell}{4}-\frac{1}{2}+\varepsilon}Q^{\frac{3}{2}}\Delta^{\frac{1}{2}}+
X^{\frac{\ell+1}{2}+\varepsilon}Q^{-1}.
\eea
Note that the first term is smaller that the third term for $\ell\geq 2$ and $Q<X$.
Take
$
Q=X^{\frac{\ell+4}{10}}\Delta^{-\frac{1}{5}}.
$
Then
\bea
\mathcal{S}_h(X)\ll_{f,\ell,\varepsilon} X^{\frac{4\ell+1}{10}+\varepsilon}\Delta^{\frac{1}{5}}.
\eea
By (3.1), (3.2) and (3.25),
\bea
\sum_{X/2<n\leq X}\lambda_f(n+h)r_{\ell}(n)\ll_{f,\ell,\varepsilon}
X^{\frac{\ell}{2}+\theta+\varepsilon}\Delta^{-1}
+ X^{\frac{4\ell+1}{10}+\varepsilon}\Delta^{\frac{1}{5}}.
\eea
By taking $\Delta=X^{\frac{\ell-1}{12}+\frac{5\theta}{6}}$ in (3.26) we obtain
\bna
\sum_{X/2<n\leq X}\lambda_f(n+h)r_{\ell}(n)\ll_{f,\ell,\varepsilon}
X^{\frac{\ell}{2}-\frac{\ell-1-2\theta}{12}+\varepsilon}.
\ena
This proves Theorem 1.

\medskip

\section{Proof of Theorem 2}
\setcounter{equation}{0}
\medskip

Denote the sum in Theorem 2 by $\mathcal{S}_h^*(X)$, By the Hardy-Littlewood-Kloosterman circle method
(see for example, \cite{I2}, Section 11.4),
we have
\bna
\mathcal{S}_h^*(X)=\int_0^1\mathscr{F}^2(\alpha)\mathscr{G}(\alpha)\mathrm{d}\alpha,
\ena
where
\bna
\mathscr{F}(\alpha)=\sum_{|m|\leq \sqrt{X}}e(\alpha m^2)
\ena
and
\bea
\mathscr{G}(\alpha)=\sum_{n\geq 1}\lambda_f(n+h)e(-\alpha n)\phi\left(\frac{n}{X}\right).
\eea
Note that $\mathscr{F}^2(\alpha)\mathscr{G}(\alpha)$ is a periodic function of period 1. We have
\bna
\mathcal{S}_h^*(X)=\int_{-1/(Q+1)}^{1-1/(Q+1)}\mathscr{F}^2(\alpha)\mathscr{G}(\alpha)\mathrm{d}\alpha,
\ena
where $Q=[5\sqrt{X}]$.
Dissecting the unit interval with Farey's points of order $Q$, we have
\bna
\mathcal{S}_h^*(X)=\sum_{q\leq Q}\sideset{}{^*}\sum_{a=1}^q
\int\limits_{\mathcal{M}(a,q)}\mathscr{F}^2\left(\frac{a}{q}+\beta\right)\mathscr{G}\left(\frac{a}{q}+\beta\right)\mathrm{d}\beta,
\ena
where the $*$ denotes that the sum is restricted by the condition $(a,q)=1$,
$
\mathcal{M}(a,q)=\left[-\frac{1}{q(q+q')},\frac{1}{q(q+q'')}\right],
$
and $\frac{a'}{q'}$, $\frac{a}{q}$ and $\frac{a''}{q''}$ are consecutive Farey fractions
which are determined
by the conditions
\[
Q<q+q',q+q''\leq q+Q, \quad aq'\equiv 1(\bmod q), \quad aq''\equiv -1(\bmod q).
\]
Exchanging the order of the summation over $a$ and the integration over $\beta$ as in
Heath-Brown \cite{HB} (see Lemma 7), we have
\bea
\mathcal{S}_h^*(X)&=&\sum_{q\leq Q}\int\limits_{|\beta|\leq \frac{1}{q Q}}\sum_{v\bmod q}
\varrho(v,q,\beta)\sum_{a=1 \atop (a,q)=1}^qe\left(-\frac{\overline{a}v}{q}\right)
\mathscr{F}^2\left(\frac{a}{q}+\beta\right)\mathscr{G}\left(\frac{a}{q}+\beta\right)\mathrm{d}\beta,
\eea
where $\varrho(v,q,\beta)$ satisfies
\bea
\varrho(v,q,\beta)\ll \frac{1}{1+|v|}
\eea
Here the implied constant is absolute.

For an asymptotic formula for $\mathscr{F}\left(\frac{a}{q}+\beta\right)$,
we quote the following result (see \cite{V}, Theorem 4.1).

\noindent {\bf Lemma 6.} {\it Suppose $(a,q)=1$, $q\leq Q$ and
$|\beta|\leq \frac{1}{qQ}$. Then we have
\bea
\mathscr{F}\left(\frac{a}{q}+\beta\right)=\frac{2G(a,0;q)}{q}\Phi_0(\beta)
+\sum_{-\frac{3q}{2}<b\leq \frac{3q}{2}}G(a,b;q)\Phi(b,q,\beta),
\eea
where $G(a,b;q)$ is the Gauss sum
\bea
G(a,b;q)=\sum\limits_{x\bmod q}e\left(\frac{ax^2+bx}{q}\right),
\eea
$\Phi_0(\beta)$ is the integral
\bea
\Phi_0(\beta)=\int_0^{\sqrt{X}}e(\beta x^2)\mathrm{d}x,
\eea
and $\Phi(b,q,\beta)$ satisfies
\bea
\sum_{-\frac{3q}{2}<b\leq \frac{3q}{2}}|\Phi(b,q,\beta)|\ll \log(q+2).
\eea
}

By (4.1) and Lemma 1, we have
\bea
\mathscr{G}\left(\frac{a}{q}+\beta\right)
&=&e\left(\frac{ha}{q}+h \beta\right)\sum_{m\geq 1}\lambda_f(m)e\left(-\frac{a m}{q}\right)
\phi\left(\frac{m-h}{X}\right)e(-\beta m)\nonumber\\
&=&e\left(\frac{ha}{q}+h \beta\right)\frac{1}{q}\sum_{m\geq 1}
\lambda_f(m)e\left(\frac{\overline{a}m}{q}\right)
\mathcal {V}_{\beta}^+\left(\frac{m}{q^2}\right)\nonumber\\
&&+\omega_f^*e\left(\frac{ha}{q}+h \beta\right)
\frac{1}{q}\sum_{m\geq 1}\lambda_f(m)e\left(-\frac{\overline{a}m}{q}\right)
\mathcal {V}_{\beta}^-\left(\frac{m}{q^2}\right),
\eea
where $\omega_f^*$ is a constant depending only on $f$ and
\bea
\mathcal {V}_{\beta}^{\pm}\left(y\right)=\int_{0}^{\infty}
\phi\left(\frac{x-h}{X}\right)e(-\beta x)\mathcal {H}_f^{\pm}(4\pi\sqrt{xy})\mathrm{d}x
\eea
with $\mathcal {H}_f^{\pm}(x)$ defined in (2.1)-(2.3).
Denote
\bea
\mathscr{G}^{\pm}\left(\frac{a}{q}+\beta\right)
=e\left(\frac{ha}{q}+h \beta\right)
\frac{1}{q}\sum_{m\geq 1}\lambda_f(m)e\left(\pm\frac{\overline{a}m}{q}\right)
\mathcal {V}_{\beta}^{\pm}\left(\frac{m}{q^2}\right),
\eea
Then by (4.4) and (4.10), we have
\bea
\sideset{}{^*}\sum_{a=1}^qe\left(-\frac{\overline{a}v}{q}\right)
\mathscr{F}^2\left(\frac{a}{q}+\beta\right)\mathscr{G}^{\pm}\left(\frac{a}{q}+\beta\right)
=\sum_{j=1}^3\mathscr{D}_j^{\pm}(v,q,\beta),
\eea
where
\bna
\mathscr{D}_1^{\pm}(v,q,\beta)&=&\frac{4e(h \beta)\Phi_0(\beta)^2}{q^3}
\sum_{m\geq 1}\lambda_f(m)\mathcal {V}_{\beta}^{\pm}\left(\frac{m}{q^2}\right)
\mathscr{C}(0,0,h,\pm m-v;q)\\
\mathscr{D}_2^{\pm}(v,q,\beta)&=&\frac{4e(h \beta)\Phi_0(\beta)}{q^2}
\sum_{m\geq 1}\lambda_f(m)\mathcal {V}_{\beta}^{\pm}\left(\frac{m}{q^2}\right)
\sum_{-\frac{3q}{2}<b\leq \frac{3q}{2}}\Phi(b,q,\beta)
\mathscr{C}(0,b,h,\pm m-v;q),\\
\mathscr{D}_3^{\pm}(v,q,\beta)&=&\frac{e(h \beta)}{q}
\sum_{m\geq 1}\lambda_f(m)\mathcal {V}_{\beta}^{\pm}\left(\frac{m}{q^2}\right)
\sum_{-\frac{3q}{2}<b_1\leq \frac{3q}{2} \atop -\frac{3q}{2}<b_2\leq \frac{3q}{2}}
 \Phi(b_1,q,\beta)\Phi(b_2,q,\beta)
\mathscr{C}(b_1,b_2,h,\pm m-v;q)
\ena
with
\bna
\mathscr{C}(b_1,b_2,h,u;q)=\sideset{}{^*}\sum_{a=1}^qG(a,b_1;q)G(a,b_2;q)e\left(\frac{a h+\overline{a}u}{q}\right),
\quad u\in \mathbb{Z}.
\ena
The following proposition will be proved in the next section.

\noindent {\bf Proposition 1.} {\it We have
\bna
\sum_{m\geq 1}|\lambda_f(m)|\left|\mathcal {V}_{\beta}^{\pm}\left(\frac{m}{q^2}\right)\right|
\ll_{f}q^2(1+|\beta|X)+\Delta q^{\frac{5}{2}}X^{-\frac{1}{4}}.
\ena}

\medskip

For the exponential sum $\mathscr{C}(b_1,b_2,h,u;q)$, we have the following estimate.

\noindent {\bf Proposition 2.} {\it Let $q=q_1q_2$, $(2q_1,q_2)=1$, $4q_1$ squarefull and $q_2$ squarefree. We have
\bna
\mathscr{C}(b_1,b_2,h,u;q)
\ll  q_1^2q_2^{\frac{3}{2}}(h,q_2)^{\frac{1}{2}}.
\ena
}

\noindent{\bf Proof.} Let $q=q_1q_2$, where $(2q_1,q_2)=1$,
$4q_1$ squarefull, $q_2$ square-free. Then we have
\bea
\mathscr{C}(b_1,b_2,h,u;q)
=\mathscr{C}(b_1,b_2,\overline{q}_2^2h,u;q_1)\mathscr{C}(b_1,b_2,\overline{q}_1^2h,u;q_2).
\eea

By Lemma 5.4.5 in \cite{Huxley}, we have $G(a,b;q)\ll \sqrt{q}$. Thus
\bea
\mathscr{C}(b_1,b_2,\overline{q}_2^2h,u;q_1)\ll q_1^2.
\eea

To estimate $\mathscr{C}(b_1,b_2,\overline{q}_1^2h,u;q_2)$, we factor
$q_2$ as $q_2=p_1p_2\cdot\cdot\cdot p_s$, $p_i$ prime. Then
\bna
\mathscr{C}(b_1,b_2,\overline{q}_1^2h,u;q_2)=\prod_{i=1}^s\mathscr{C}(b_1,b_2,\overline{p_i'}^2\overline{q}_1^2h,u;p_i),
\ena
where $p_i'=q_2/p_i$. Thus we consider the exponential sum
\bna
\mathscr{C}(b_1,b_2,rh,u;p)=\sideset{}{^*}\sum_{a=1}^pG(a,b_1;p)G(a,b_2;p)e\left(\frac{a rh+\overline{a}u}{p}\right),
\quad (p,2r)=1,\quad  p\,\mbox{prime}.
\ena
By Lemma 5.4.5 in \cite{Huxley}, we have, for $(2a,q)=1$,
\bna
G(a,b;q)=
e\left(-\frac{\bar{4}\bar{a}b^2}{q}\right)G(a,0;q)
=e\left(-\frac{\bar{4}\bar{a}b^2}{q}\right)\left(\frac{a}{q}\right)\epsilon_q\sqrt{q},
\ena
where $\epsilon_q$ is defined in (2.7).
Thus by Weil's bound for Kloosterman sum, we have
\bna
\mathscr{C}(b_1,b_2,rh,u;p)&=&\sideset{}{^*}\sum_{a=1}^p
e\left(-\frac{\bar{4}\bar{a}b_1^2}{p}\right)\left(\frac{a}{p}\right)\epsilon_p\sqrt{p}
e\left(-\frac{\bar{4}\bar{a}b_2^2}{p}\right)\left(\frac{a}{p}\right)\epsilon_p\sqrt{p}
e\left(\frac{a rh+\overline{a}u}{p}\right)\\
&=&p\epsilon_p^2\sideset{}{^*}\sum_{a=1}^p
e\left(\frac{a rh+\overline{a}(u-\bar{4}b_1^2-\bar{4}b_2^2)}{p}\right)\\
&\ll&p^{\frac{3}{2}}(rh,u-\bar{4}b_1^2-\bar{4}b_2^2,p)^{\frac{1}{2}}\\
&\ll&p^{\frac{3}{2}}(h,p)^{\frac{1}{2}}.
\ena
It follows that
\bea
\mathscr{C}(b_1,b_2,\overline{q}_1^2h,u;q_2)\ll q_2^{\frac{3}{2}}(h,q_2)^{\frac{1}{2}}.
\eea
By (4.12)-(4.14), Proposition 2 follows. \hfill $\Box$

\medskip

By the second derivative test for exponential integrals and the trivial estimate,
$\Phi_0(\beta)$ in (4.6) is bounded by
\bea
\Phi_0(\beta)\ll \left(\frac{X}{1+|\beta|X}\right)^{\frac{1}{2}}.
\eea
By (4.3), (4.15) and Propositions 1-2, we have, for $q=q_1q_2$, $(2q_1,q_2)=1$,
$4q_1$ squarefull and $q_2$ square-free,
\bna
&&\int\limits_{|\beta|\leq \frac{1}{q Q}}\sum_{v\bmod q}
\varrho(v,q,\beta)\mathscr{D}_1^{\pm}(v,q,\beta)\mathrm{d}\beta\nonumber\\
&\ll_f& (\log Q)
 \frac{1}{q_1^3q_2^3}\int\limits_{|\beta|\leq \frac{1}{q_1q_2 Q}}\frac{X}{1+|\beta|X}
\left(q_1^2q_2^2(1+|\beta|X)+\Delta q_1^{\frac{5}{2}}q_2^{\frac{5}{2}}X^{-\frac{1}{4}}\right)
q_1^2q_2^{\frac{3}{2}}(h,q_2)^{\frac{1}{2}}\mathrm{d}\beta\nonumber\\
&\ll_f&\left(\log X\right)^2
 \left(\frac{Xq_1^2q_2^2}{q_1q_2Q}+\Delta q_1^{\frac{5}{2}}q_2^{\frac{5}{2}}X^{-\frac{1}{4}} \right)
 q_1^{-1}q_2^{-\frac{3}{2}}(h,q_2)^{\frac{1}{2}},
\ena
and
\bea
&&\sum_{q\leq Q}\int\limits_{|\beta|\leq \frac{1}{q Q}}\sum_{v\bmod q}
\varrho(v,q,\beta)\mathscr{D}_1^{\pm}(v,q,\beta)\mathrm{d}\beta\nonumber\\
&\ll_{f,\varepsilon}&\frac{X^{1+\varepsilon}}{Q}\sum_{q_2\leq Q}(h,q_2)^{\frac{1}{2}}q_2^{-\frac{1}{2}}
\sum_{q_1\leq Q/q_2 \atop 4q_1 \mathrm{squarefull}}1
+\Delta X^{-\frac{1}{4}+\varepsilon}\sum_{q_2\leq Q}(h,q_2)^{\frac{1}{2}}q_2
\sum_{q_1\leq Q/q_2 \atop 4q_1 \mathrm{squarefull}}q_1^{\frac{3}{2}}\nonumber\\
&\ll_{f,\varepsilon}&\frac{X^{1+\varepsilon}}{Q}
\sum_{q_2\leq Q}(h,q_2)^{\frac{1}{2}}q_2^{-\frac{1}{2}}
\left(\frac{Q}{q_2}\right)^{\frac{1}{2}}
+\Delta X^{-\frac{1}{4}+\varepsilon}
\sum_{q_2\leq Q}(h,q_2)^{\frac{1}{2}}q_2\left(\frac{Q}{q_2}\right)^2\nonumber\\
&\ll_{f,\varepsilon}&\frac{X^{1+\varepsilon}}{\sqrt{Q}}\sum_{d_1|h}\frac{1}{\sqrt{d_1}}\sum_{d_2\leq Q/d_1}\frac{1}{d_2}
+\Delta X^{-\frac{1}{4}+\varepsilon}Q^2\sum_{d_1|h}\frac{1}{\sqrt{d_1}}\sum_{d_2\leq Q/d_1}\frac{1}{d_2}\nonumber\\
&\ll_{f,\varepsilon}&\frac{h^{\varepsilon}X^{1+\varepsilon}}{\sqrt{Q}}+
h^{\varepsilon}\Delta X^{-\frac{1}{4}+\varepsilon}Q^2\nonumber\\
&\ll_{f,\varepsilon}&\Delta X^{\frac{3}{4}+\varepsilon},
\eea
uniformly for $1\leq h\leq X$.

Similarly, by (4.3), (4.7), (4.15) and Propositions 1-2, we have
\bna
&&\int\limits_{|\beta|\leq \frac{1}{q Q}}\sum_{v\bmod q}
\varrho(v,q,\beta)\mathscr{D}_2^{\pm}(v,q,\beta)\mathrm{d}\beta\nonumber\\
&\ll_f&
\frac{(\log Q)^2}{q_1^2q_2^2}\int\limits_{|\beta|\leq \frac{1}{q_1q_2 Q}}\left(\frac{X}{1+|\beta|X}\right)^{\frac{1}{2}}
\left(q_1^2q_2^2(1+|\beta|X)+\Delta q_1^{\frac{5}{2}}q_2^{\frac{5}{2}}X^{-\frac{1}{4}}\right)
q_1^2q_2^{\frac{3}{2}}(h,q_2)^{\frac{1}{2}}\mathrm{d}\beta\nonumber\\
&\ll_{f}&(\log X)^2
\left(\frac{X^{\frac{1}{2}}q_1^2q_2^2}{q_1q_2Q}+\frac{Xq_1^2q_2^2}{(q_1q_2Q)^{\frac{3}{2}}}
+\frac{\Delta q_1^{\frac{5}{2}}q_2^{\frac{5}{2}}X^{-\frac{1}{4}}}{\sqrt{q_1q_2Q}}\right)
q_2^{-\frac{1}{2}}(h,q_2)^{\frac{1}{2}},
\ena
and
\bea
&&\sum_{q\leq Q}\int\limits_{|\beta|\leq \frac{1}{q Q}}\sum_{v\bmod q}
\varrho(v,q,\beta)\mathscr{D}_2^{\pm}(v,q,\beta)\mathrm{d}\beta\nonumber\\
&\ll_{f,\varepsilon}&\frac{X^{\frac{1}{2}+\varepsilon}}{Q}
\sum_{q_2\leq Q}(h,q_2)^{\frac{1}{2}}q_2^{\frac{1}{2}}
\sum_{q_1\leq Q/q_2 \atop 4q_1 \mathrm{squarefull}}q_1
+\frac{X^{1+\varepsilon}}{Q^{\frac{3}{2}}}\sum_{q_2\leq Q}(h,q_2)^{\frac{1}{2}}
\sum_{q_1\leq Q/q_2 \atop 4q_1 \mathrm{squarefull}}q_1^{\frac{1}{2}}\nonumber\\
&&+\Delta X^{-\frac{1}{4}+\varepsilon}Q^{-\frac{1}{2}}\sum_{q_2\leq Q}(h,q_2)^{\frac{1}{2}}q_2^{\frac{3}{2}}
\sum_{q_1\leq Q/q_2 \atop 4q_1 \mathrm{squarefull}}q_1^2\nonumber\\
&\ll_{f,\varepsilon}&\frac{X^{\frac{1}{2}+\varepsilon}}{Q}
\sum_{q_2\leq Q}(h,q_2)^{\frac{1}{2}}q_2^{\frac{1}{2}}
\left(\frac{Q}{q_2}\right)^{\frac{3}{2}}
+\frac{X^{1+\varepsilon}}{Q^{\frac{3}{2}}}
\sum_{q_2\leq Q}(h,q_2)^{\frac{1}{2}}\left(\frac{Q}{q_2}\right)\nonumber\\
&&+\Delta X^{-\frac{1}{4}+\varepsilon}Q^{-\frac{1}{2}}
\sum_{q_2\leq Q}(h,q_2)^{\frac{1}{2}}q_2^{\frac{3}{2}}
\left(\frac{Q}{q_2}\right)^{\frac{5}{2}}\nonumber\\
&\ll_{f,\varepsilon}&X^{\frac{1}{2}+\varepsilon}Q^{\frac{1}{2}}
\sum_{d_1|h}d_1^{-\frac{1}{2}}\sum_{d_2\leq Q/d_1}\frac{1}{d_2}
+\frac{X^{1+\varepsilon}}{Q^{\frac{1}{2}}}
\sum_{d_1|h}d_1^{-\frac{1}{2}}\sum_{d_2\leq Q/d_1}\frac{1}{d_2}\nonumber\\
&&+\Delta X^{-\frac{1}{4}+\varepsilon}Q^2\sum_{d_1|h}\frac{1}{\sqrt{d_1}}\sum_{d_2\leq Q/d_1}\frac{1}{d_2}\nonumber\\
&\ll_{f,\varepsilon}&h^{\varepsilon}X^{\frac{1}{2}+\varepsilon}Q^{\frac{1}{2}}
+\frac{h^{\varepsilon}X^{1+\varepsilon}}{Q^{\frac{1}{2}}}
+h^{\varepsilon}\Delta X^{-\frac{1}{4}+\varepsilon}Q^2\nonumber\\
&\ll_{f,\varepsilon}&\Delta X^{\frac{3}{4}+\varepsilon},
\eea
uniformly for $1\leq h\leq X$.

Further, by (4.3), (4.7), (4.15) and Propositions 1-2, we have
\bna
&&\int\limits_{|\beta|\leq \frac{1}{q Q}}\sum_{v\bmod q}
\varrho(v,q,\beta)\mathscr{D}_3^{\pm}(v,q,\beta)\mathrm{d}\beta\nonumber\\
&\ll_f&
 \frac{(\log Q)^3}{q_1q_2}\int\limits_{|\beta|\leq \frac{1}{q_1q_2 Q}}
\left(q_1^2q_2^2(1+|\beta|X)+\Delta q_1^{\frac{5}{2}}q_2^{\frac{5}{2}}X^{-\frac{1}{4}}\right)
q_1^2q_2^{\frac{3}{2}}(h,q_2)^{\frac{1}{2}}\mathrm{d}\beta\nonumber\\
&\ll_{f}&(\log X)^3
\left(\frac{q_1^2q_2^2}{q_1q_2Q}+\frac{Xq_1^2q_2^2}{(q_1q_2Q)^2}
+\frac{\Delta q_1^{\frac{5}{2}}q_2^{\frac{5}{2}}X^{-\frac{1}{4}}}{q_1q_2Q}\right)
q_1q_2^{\frac{1}{2}}(h,q_2)^{\frac{1}{2}},
\ena
and
\bea
&&\sum_{q\leq Q}\int\limits_{|\beta|\leq \frac{1}{q Q}}\sum_{v\bmod q}
\varrho(v,q,\beta)\mathscr{D}_3^{\pm}(v,q,\beta)\mathrm{d}\beta\nonumber\\
&\ll_{f,\varepsilon}&\frac{X^{\varepsilon}}{Q}
\sum_{q_2\leq Q}(h,q_2)^{\frac{1}{2}}q_2^{\frac{3}{2}}
\sum_{q_1\leq Q/q_2 \atop 4q_1 \mathrm{squarefull}}q_1^2
+\frac{X^{1+\varepsilon}}{Q^2}\sum_{q_2\leq Q}(h,q_2)^{\frac{1}{2}}q_2^{\frac{1}{2}}
\sum_{q_1\leq Q/q_2 \atop 4q_1 \mathrm{squarefull}}q_1\nonumber\\
&&+\Delta X^{-\frac{1}{4}+\varepsilon}Q^{-1}\sum_{q_2\leq Q}(h,q_2)^{\frac{1}{2}}q_2^2
\sum_{q_1\leq Q/q_2 \atop 4q_1 \mathrm{squarefull}}q_1^{\frac{5}{2}}\nonumber\\
&\ll_{f,\varepsilon}&\frac{X^{\varepsilon}}{Q}
\sum_{q_2\leq Q}(h,q_2)^{\frac{1}{2}}q_2^{\frac{3}{2}}
\left(\frac{Q}{q_2}\right)^{\frac{5}{2}}
+\frac{X^{1+\varepsilon}}{Q^2}\sum_{q_2\leq Q}
(h,q_2)^{\frac{1}{2}}q_2^{\frac{1}{2}}\left(\frac{Q}{q_2}\right)^{\frac{3}{2}}\nonumber\\
&&+\Delta X^{-\frac{1}{4}+\varepsilon}Q^{-1}\sum_{q_2\leq Q}(h,q_2)^{\frac{1}{2}}q_2^2
\left(\frac{Q}{q_2}\right)^3\nonumber\\
&\ll_{f,\varepsilon}&X^{\varepsilon}Q^{\frac{3}{2}}
\sum_{d_1|h}d_1^{-\frac{1}{2}}\sum_{d_2\leq Q/d_1}\frac{1}{d_2}
+\frac{X^{1+\varepsilon}}{Q^{\frac{1}{2}}}\sum_{d_1|h}d_1^{-\frac{1}{2}}\sum_{d_2\leq Q/d_1}\frac{1}{d_2}
+\Delta X^{-\frac{1}{4}+\varepsilon}Q^2\sum_{d_1|h}\frac{1}{\sqrt{d_1}}\sum_{d_2\leq Q/d_1}\frac{1}{d_2}\nonumber\\
&\ll_{f,\varepsilon}&h^{\varepsilon}X^{\varepsilon}Q^{\frac{3}{2}}
+\frac{h^{\varepsilon}X^{1+\varepsilon}}{Q^{\frac{1}{2}}}
+h^{\varepsilon}\Delta X^{-\frac{1}{4}+\varepsilon}Q^2\nonumber\\
&\ll_{f,\varepsilon}&\Delta X^{\frac{3}{4}+\varepsilon},
\eea
uniformly for $1\leq h\leq X$.
By (4.2), (4.8), (4.11) and (4.16)-(4.18), Theorem 2 follows.

\medskip

\section{Proof of Proposition 1}
\setcounter{equation}{0}
\medskip

Recall $\mathcal {V}_{\beta}^{\pm}(y)$ in (4.9) which we relabel as
\bea
\mathcal {V}_{\beta}^{\pm}\left(y\right)=\int_{0}^{\infty}
\phi\left(\frac{x-h}{X}\right)e(-\beta x)\mathcal {H}_f^{\pm}(4\pi\sqrt{xy})\mathrm{d}x,
\eea
where $\mathcal {H}_f^{\pm}(x)$ defined in (2.1)-(2.3).
Note that for $y=m/q^2$, $q\leq Q=[5\sqrt{X}]$ and $x\gg X$, we have $4\pi\sqrt{xy}>1$.
By Lemma 5, we have
\bea
\mathcal {H}_f^-\left(4\pi\sqrt{xy}\right)\ll_f x^{-\frac{5}{4}}y^{-\frac{5}{4}},
\eea
and
\bea
\mathcal {H}_f^+\left(4\pi\sqrt{xy}\right)=(x y)^{-\frac{1}{4}}\sum_{\pm}c_1^{\pm}e\left(\pm 2\sqrt{x y}\right)
+(x y)^{-\frac{3}{4}}\sum_{\pm}c_2^{\pm}e\left(\pm 2\sqrt{x y}\right)
+O_f\left(x^{-\frac{5}{4}}y^{-\frac{5}{4}}\right),
\eea
where $c_1^{\pm}$ and $c_2^{\pm}$ are absolute constants depending only on the weight $\kappa$ or the
spectral parameter $\mu$.
Plugging (5.2) and (5.3) into (5.1), we have
\bea
\mathcal {V}_{\beta}^{-}\left(y\right)\ll y^{-\frac{5}{4}}\int_0^{\infty}
x^{-\frac{5}{4}}\phi\left(\frac{x-h}{X}\right)\mathrm{d}x
\ll y^{-\frac{5}{4}}X^{-\frac{1}{4}}
\eea
and
\bna
\mathcal {V}_{\beta}^{+}(y)
&=& y^{-\frac{1}{4}}\sum_{\pm}c_1^{\pm}\int_0^{\infty}x^{-\frac{1}{4}}
\phi\left(\frac{x-h}{X}\right)e\left(-\beta x\pm 2\sqrt{xy}\right)\mathrm{d}x\\
&&+y^{-\frac{3}{4}}\sum_{\pm}c_2^{\pm}\int_0^{\infty}x^{-\frac{3}{4}}
\phi\left(\frac{x-h}{X}\right)e\left(-\beta x\pm 2\sqrt{xy}\right)\mathrm{d}x
+O_f\left(y^{-\frac{5}{4}}X^{-\frac{1}{4}}\right).
\ena
Changing variable $x\rightarrow Xt^2$, we have
\bea
\mathcal {V}_{\beta}^{+}\left(y\right)
&=&2y^{-\frac{1}{4}}X^{\frac{3}{4}}\sum_{\pm}c_1^{\pm}
\int_0^{\infty}t^{\frac{1}{2}}
\phi\left(t^2-\frac{h}{X}\right)e\left(-\beta Xt^2\pm 2\sqrt{yX}t\right)\mathrm{d}t\nonumber\\
&&+2y^{-\frac{3}{4}}X^{\frac{1}{4}}\sum_{\pm}c_2^{\pm}
\int_0^{\infty}t^{-\frac{1}{2}}
\phi\left(t^2-\frac{h}{X}\right)e\left(-\beta Xt^2\pm 2\sqrt{yX}t\right)\mathrm{d}x
+O_f\left(y^{-\frac{5}{4}}X^{-\frac{1}{4}}\right).\nonumber
\\
\eea
Denote $\rho(t)=-\beta Xt^2\pm 2\sqrt{yX}t$.
If $y>4|\beta|^2X$, we have
$|\rho'(t)|=|-2\beta Xt\pm 2\sqrt{yX}|\geq 2\sqrt{yX}-2|\beta|Xt\geq 2\sqrt{yX}-2|\beta|X\geq \sqrt{yX}$.
By partial integration twice, we have
\bna
\int_0^{\infty}t^{\pm\frac{1}{2}}
\phi\left(t^2-\frac{h}{X}\right)e\left(-\beta Xt^2\pm 2\sqrt{yX}t\right)\mathrm{d}t\ll \Delta y^{-1}X^{-1},
\ena
and by (5.5),
\bea
\mathcal {V}_{\beta}^{+}(y)\ll_{f}\Delta y^{-\frac{5}{4}}X^{-\frac{1}{4}}.
\eea

If $y\leq 4|\beta|^2X$, we have $|\rho''(t)|=|-2\beta X|\gg 1$. By the second derivative test,
\bna
\int_0^{\infty}t^{\pm\frac{1}{2}}
\phi\left(t^2-\frac{h}{X}\right)e\left(-\beta Xt^2\pm 2\sqrt{yX}t\right)\mathrm{d}t\ll \frac{1}{\sqrt{|\beta|X}},
\ena
and by (5.5),
\bea
\mathcal {V}_{\beta}^{+}(y)\ll_{f} y^{-\frac{1}{4}}X^{\frac{1}{4}}|\beta|^{-\frac{1}{2}}+
y^{-\frac{3}{4}}X^{-\frac{1}{4}}|\beta|^{-\frac{1}{2}}
+y^{-\frac{5}{4}}X^{-\frac{1}{4}}.
\eea
By (3.17), (5.4), (5.6) and (5.7), we have
\bna
&&\sum_{m\geq 1}|\lambda_f(m)|\left|\mathcal {V}_{\beta}^{\pm}\left(\frac{m}{q^2}\right)\right|\\
&\ll_f&
\sum_{m\leq 4|\beta|^2Xq^2}|\lambda_f(m)|
\left(\left(\frac{m}{q^2}\right)^{-\frac{1}{4}}X^{\frac{1}{4}}|\beta|^{-\frac{1}{2}}+
\left(\frac{m}{q^2}\right)^{-\frac{3}{4}}X^{-\frac{1}{4}}|\beta|^{-\frac{1}{2}}
+\left(\frac{m}{q^2}\right)^{-\frac{5}{4}}X^{-\frac{1}{4}}\right)\\
&&+
\Delta\sum_{m> 4|\beta|^2Xq^2}|\lambda_f(m)|\left(\frac{m}{q^2}\right)^{-\frac{5}{4}}X^{-\frac{1}{4}}\\
&\ll_f&q^{\frac{1}{2}}|\beta|^{-\frac{1}{2}}X^{\frac{1}{4}}\sum_{m\leq 4|\beta|^2Xq^2}|\lambda_f(m)|m^{-\frac{1}{4}}
+q^{\frac{3}{2}}|\beta|^{-\frac{1}{2}}X^{-\frac{1}{4}}\sum_{m\leq 4|\beta|^2Xq^2}
|\lambda_f(m)|m^{-\frac{3}{4}}\\
&&+q^{\frac{5}{2}}X^{-\frac{1}{4}}\sum_{m\leq 4|\beta|^2Xq^2}|\lambda_f(m)|m^{-\frac{5}{4}}
+\Delta q^{\frac{5}{2}}X^{-\frac{1}{4}}
\\
&\ll_{f}&q^2|\beta|X+q^2+\Delta q^{\frac{5}{2}}X^{-\frac{1}{4}}.
\ena
This proves Proposition 1.

\bigskip

\noindent
{\sc Acknowledgements.}
The author is very grateful to the referee
for detailed comments and valuable suggestions which bring many improvements
on the original draft.
This work is supported by
National Natural Science Foundation of China (11101239),
Young Scholars Program of Shandong University, Weihai (2015WHWLJH04) and
Natural Science Foundation of Shandong Province (ZR2016AQ15).

\end{document}